\documentclass[12pt]{amsart}

\usepackage[alphabetic]{amsrefs}

\usepackage{amsthm,amscd}
\theoremstyle{plain}
\newtheorem{thm}{Theorem}[section]
\newtheorem*{thm*}{Theorem}

\newtheorem{prop}[thm]{Proposition}
\newtheorem*{prop*}{Proposition}

\theoremstyle{definition}

\theoremstyle{remark}
\newtheorem*{rem}{Remark}

\newcommand{\grd}[2]{\mathfrak{g}^{#1}_{#2}}

\newcommand{\Pic}{\operatorname{Pic}}

\newcommand{\mult}{\operatorname{mult}}

\newcommand{\C}{\mathbf{C}} 
\newcommand{\OO}{\mathcal{O}} 
\newcommand{\Z}{\mathbf{Z}} 
\newcommand{\A}{\mathbf{A}}

\newcommand{\PR}{\mathbf{P}}
\newcommand{\M}{\mathcal{M}}

\def\topbotatom#1{\hbox{\hbox to 0pt{$#1\bot$\hss}$#1\top$}}

\usepackage[left=2cm,top=2cm,right=2cm]{geometry}                
\usepackage{graphicx}
\usepackage{amssymb,amsmath, gensymb}
\usepackage{epstopdf}
\usepackage{hyperref, enumerate}

\newcommand{\fun}[3]{#1\!\!:#2\rightarrow #3}

\title[Irreducibility and stable rationality of Weierstrass points in low genus]{Irreducibility and stable rationality of the loci of Weierstrass points on curves of genus at most six}
\author{Evan M. Bullock}

\begin{document}
\maketitle
\begin{abstract} Given a numerical semigroup $H\subseteq (\Z_{\geq 0},+)$, we consider the locus $\M_{g,1}^H$ of smooth curves of genus $g$ with a marked Weierstrass point of semigroup $H$.  We show that for all semigroups $H$ of genus $g\leq 6$ the locus $\M_{g,1}^H$ is irreducible and that for all but possibly two such semigroups it is stably rational.
\end{abstract}

\section{Introduction}

By a \emph{curve}, we will mean a projective curve over the complex numbers.  Given a smooth connected curve $C$ of genus $g>1$ and a point $p\in C$, the \emph{Weierstrass semigroup} of $C$ at $p$ is the set
\[H=\{n\in \Z_{\geq 0}\!:~\text{some meromorphic function $f$ has a pole of order $n$ at $p$ and no other poles}\}.\]
It follows from the Riemann-Roch theorem that if $H\subseteq\Z_{\geq 0}$ is the Weierstrass semigroup of $C$ at $p$, then $H$ consists of all but $g$ of the positive integers, called the \emph{Weierstrass gaps}.  If $H\subseteq (\Z_{\geq 0},+)$ is any subsemigroup with $|\Z_{\geq 0}\backslash H|=g$, we say that $H$ is a \emph{numerical semigroup of genus $g$}.

For a general point $p\in C$, the Weierstrass semigroup is $\langle g+1,g+2,g+3,\ldots\rangle$ and the gaps are $1,2,\ldots,g$.  A \emph{Weierstrass point} on $C$ is any point with a different semigroup.  The pairs $(C,p)$ where $p$ is a Weierstrass point of a smooth irreducible curve $C$ form a divisor in the coarse moduli space $\M_{g,1}$ of pointed curves.  The subset $\M_{g,1}^H\subseteq\M_{g,1}$ consisting of $(C,p)$ where $p$ has semigroup $H$ is locally closed.

An irreducible variety $V$ of dimension $n$ is \emph{stably rational} if $V\times \PR^m$ is birational to $\PR^{n+m}$ for sufficiently large $m$; this condition is strictly stronger than unirationality, but strictly weaker than rationality.  
For example, stably rational varieties that are not rational are constructed in \cite{stabledoesnotimplyrational} and in \cite{unirationalnotstablyrational} Artin and Mumford show that the torsion part of $H^3(X,\Z)$ is a birational invariant of a smooth projective variety $X$ and construct unirational varieties of each dimension $n\geq 3$ with non-trivial torsion in $H^3$; such varieties cannot be stably rational by the K\"{u}nneth theorem.

In this paper, we prove the following:

\begin{thm}\label{maintheorem} Let $H$ be a semigroup of genus $g\leq 6$.  Then the locus $\M^H_{g,1}$ of smooth curves of genus $g$ with a marked point of Weierstrass semigroup $H$ is irreducible.  Moreover, $\M^H_{g,1}$ is always stably rational, with the possible exceptions of the semigroup $\langle 5,7,8,9,11\rangle$ of a general Weierstrass point in genus $5$ and the semigroup $\langle 6,7,8,9,10\rangle$ of a general odd subcanonical point in genus $6$.
\end{thm}

We will prove this in Section \ref{proofofmainthm}, after recalling the equivalence between the Weierstrass semigroup of $C$ at $p$ and the vanishing and ramification sequences of the canonical series $K_C$ at $p$ and discussing previous results on the irreducibility of $\M^H_{g,1}$. 

Weierstrass points can also be defined in terms of holomorphic $1$-forms rather than meromorphic functions.  Given a point $p$ on a smooth curve $C$, the \emph{vanishing sequence} (of the canonical series $K_C$) at $p$ is the ascending sequence $0=a_0(p)<a_1(p)<\cdots<a_{g-1}(p)\leq 2g-2$ of non-negative integers such that 
\[\{a_i(p)\}=\left\{v_p(\omega)\!:~\omega\in H^0\left(C,K_C\right)\right\},\]
i.e. the $a_i(p)$ are the zero orders at $p$ of the global holomorphic $1$-forms on $C$.  It follows from the Riemann-Roch theorem (cf. \cite{ACGH} I.E) that the set of Weierstrass gaps at $p$ is $\{a_i(p)+1\}$, so a Weierstrass point is a point where the vanishing sequence is bigger than the sequence $a_i=i$ at a general point.

The \emph{ramification sequence} $0=\alpha_0\leq \alpha_1\leq\cdots\leq \alpha_{g-1}\leq g-1$ is the non-decreasing sequence defined by $\alpha_i(p)=a_i(p)-i$; it also encodes the same information as the vanishing sequence or the semigroup.  The \emph{weight} of a Weierstrass point,  \[w(p) = \sum_{i=0}^{g-1} \alpha_i(p)\]
measures how far a point is from being general.  In particular, if $H$ is a semigroup of weight $w$, then every component of $\M_{g,1}^H$ has codimension at most $w$ in $\M_{g,1}$; a component is called \emph{dimensionally proper} if the codimension is equal to the weight.

The most powerful tool for showing the existence of Weierstrass points for wide classes of numerical semigroups (i.e. the non-emptiness of $\M_{g,1}^H$) is the method of limit linear series of Eisenbud and Harris \cite{EHexist}, which has since been used to prove several slight extensions of their results (e.g. \cite{existenceweightg-1}), and was a major tool in showing that every numerical semigroup of genus at most $8$ is a Weierstrass semigroup in \cite{existenceleq8}.  This method involves attaching an elliptic tail to a Weierstrass point on a genus $g-1$ curve and smoothing to get a genus $g$ curve.  It is not typically well-suited to proving irreducibility because there is no easy way of showing that the closure of a component of $\M_{g,1}^H$ in $\overline{\M_{g,1}}$ would necessarily meet the boundary component $\Delta_{1,1}$.  

On the other hand, there are relatively few infinite classes of semigroups for which irreducibility is known (e.g. for ordinary Weierstrass points of weight $1$ this was shown in \cite{EHmono}).  Weierstrass points with first non-gap $3$ must lie on trigonal curves, and that gap sequences of such points have been completely classified, including irreducibility for each gap sequence, in \cite{trigonal}.

In \cite{Arbthesis} Arbarello proves the irreducibility of the locus $W_{n,g}$ of Weierstrass points whose first non-gap is $n$; this certainly shows that some $\M_{g,1}^H$ is irreducible, but while we should expect $H$ to be the ``smallest'' semigroup (with respect to the natural partial order on the corresponding vanishing or ramification sequences) with this property, as far as I know this has never been proven.

Similarly, in their study of the moduli space of pairs $(C,\omega)$ of a curve and a holomorphic $1$-form with a prescribed partition of zeros in \cite{KZ}, Kontsevich and Zorich in particular analyze the case of \emph{subcanonical points}, i.e. points $p\in C$ such that $K_C\cong \OO_C\left(\left(2g-2\right)p \right)$, and show that this locus has three irreducible components for all $g\geq 4$.  In \cite{mythesis}, we show that each component is a different $\M_{g,1}^H$ and determine for each genus the three semigroups, using a limit linear series argument as in \cite{EHexist}.

The general irreducibility theorems covering the most different non-trigonal semigroups, however, are still those using earlier techniques.  In \cite{Pink}, Pinkham constructed a compactification of $\M_{g,1}^H$ for each numerical semigroup $H=\langle k_1,\ldots,k_r\rangle$ by studying the weighted deformation theory of the singular monomial curve in $\A^r$ parameterized by $(t^{k_1},\ldots,t^{k_r})$.  In \cite[4.2]{Buch80}, Buchweitz noticed that applying some general unobstructedness results from deformation theory could be used to show the existence of Weierstrass points, in the cases where $H$ has two or three generators, and the case where $H$ has four generators and the ring $\C[t^{k_1},\ldots,t^{k_r}]$ is Gorenstein.
Of course if there are no obstructions, the Pinkham compactification of $\M_{g,1}^H$ is just a weighted projective space, which is certainly irreducible, so the Buchweitz result implies irreducibility as well.  

The Gorenstein property in this case was shown to have a purely combinatorial description in \cite{Gorenstein}: let $c$ be the largest gap of $H$.  We say that $H$ is \emph{symmetric} iff for all $0\leq i\leq c$, we have that $i$ is a gap if and only if $c-i$ is a non-gap.  Then $\C[t^{k_1},\ldots,t^{k_r}]$ is Gorenstein if and only if $H$ is symmetric.  We thus have the following (the $r\leq 3$ part of which was known to Nakano in \cite{Nak2}):

\begin{thm} Let $H=\langle k_1,\ldots,k_r\rangle$ be a numerical semigroup.  Suppose either that $r\leq 3$ or that $r=4$ and $H$ is symmetric.  Then $\M_{g,1}^H$ is non-empty and irreducible and its Pinkham compactification is a weighted projective space.
\end{thm}

\begin{rem} The Weierstrass points with symmetric semigroup are precisely the subcanonical points discussed in \cite{mythesis}.  This follows from \cite[Lemma]{Gorenstein} since a point is subcanonical if and only if $2g-1$ is a gap. 
\end{rem}

In \cite{Nak1} and \cite{Nak2}, Nakano and Mori used computer deformation theory calculations to prove irreducibility and rationality of $\M_{g,1}^H$ for (almost) all the four-generator semigroups $H$ of genus at most $6$.  In our proof of Theorem \ref{maintheorem} in the next section, we use other methods to study the one remaining four-generator case and all the cases where $H$ has more than four generators.

\section{Proof of Theorem \ref{maintheorem}}\label{proofofmainthm}

By the results \cite{Nak1, Nak2} mentioned above, together with the results of Ballico, Casnati, and Fontanari \cite{Mg1, M61} showing that $\M_{g,1}$ itself is rational for $g\leq 6$, the irreducibility and rationality of $\M_{g,1}^H$ is  known for all semigroups of genus $g\leq 6$ except the following nine cases:

\begin{center}
  \begin{tabular}{| l | l | l | }
    \hline
    semigroup & vanishing sequence & ramification sequence \\ \hline
    $N(5)_8=\langle 5,7,8,9,11\rangle$ & $0,1,2,3,5$ & $0,0,0,0,1$ \\ \hline\hline
    $N(6)_{5}=\langle 4,9,10,11\rangle$ & $0,1,2,4,5,6$ & $0,0,0,1,1,1$ \\ \hline
    $N(6)_{11}=\langle 5,8,9,11,12\rangle$ & $0,1,2,3,5,6$ & $0,0,0,0,1,1$ \\ \hline 
    $N(6)_{12}=\langle 5,8,9,11,13\rangle$ & $0,1,2,3,5,7$ & $0,0,0,0,1,2$ \\ \hline
    $N(6)_{18}=\langle 6,8,9,10,11,13\rangle$ & $0,1,2,3,4,6$ & $0,0,0,0,0,1$ \\ \hline
    $N(6)_{19}=\langle 6,7,9,10,11\rangle$ & $0,1,2,3,4,7$ & $0,0,0,0,0,2$ \\ \hline
    $N(6)_{20}=\langle 6,7,8,10,11\rangle$ & $0,1,2,3,4,8$ & $0,0,0,0,0,3$ \\ \hline
    $N(6)_{21}=\langle 6,7,8,9,11\rangle$ & $0,1,2,3,4,9$ & $0,0,0,0,0,4$ \\ \hline
    $N(6)_{22}=\langle 6,7,8,9,10\rangle$ & $0,1,2,3,4,10$ & $0,0,0,0,0,5$ \\ \hline
  \end{tabular}
\end{center}

The locus $\M_{6,1}^{N(6)_{5}}$ is an open subset of the locus $W_{4,6}$ of pointed curves of genus $6$ possessing a $\grd{1}{4}$ that is totally ramified at the marked point: this is the locus of pointed curves $(C,p)$ where $h^0\left(C,\OO_C(4p)\right)\geq 2$, and by Riemann-Roch this is equivalent to the condition $a_3(p)\geq 4$.  Since $0,1,2,4,5,6$ is the smallest possible vanishing sequence in genus $6$ that satisfies this condition (and Weierstrass points of any semigroup of genus at most $8$ exist by \cite{existenceleq8}), it is the vanishing sequence at a general point in $W_{4,6}$.

This locus $W_{4,g}$ has been studied for general genus $g$: it was shown to be irreducible of dimension $2g+1$ in \cite{Arbthesis} and rational in \cite{type4g}.  By the same reasoning,  $\M_{6,1}^{N(6)_{11}}$ is an open subset of $W_{5,6}$, and hence is irreducible of dimension $2g+2$ by \cite{Arbthesis}.

The semigroups $N(5)_8$ and $N(6)_{18}$ correspond to general Weierstrass points, and irreducibility is known in this case for every genus.  In fact, in \cite{EHmono} a stronger result was proved: for $H=\langle g,g+2,g+3\ldots,2g+1\rangle$ the semigroup of a general Weierstrass point of genus $g$, the monodromy of the cover $\M^H_{g,1}\to \M_g$ is the full symmetric group $S_{(g-1)g(g+1)}$.

In one other case irreducibility is known: $\M_{6,1}^{N(6)_{22}}$ is the locus of odd subcanonical points with the smallest possible semigroup, which is shown to be non-empty and irreducible for any genus in \cite{mythesis}.

To prove irreducibility and stable rationality in the remaining cases, we first recall a description of (most) genus $6$ curves:
\begin{prop}[\cite{ACGH} V.A, \cite{SB89}]\label{acghprop} Let $C\subset \PR^5$ be a canonical curve of genus $6$ which is not hyperelliptic, trigonal, bielliptic, isomorphic to a smooth plane quintic curve, or birational to a plane sextic with double point singularities, at least one of which is tacnodal (or worse).  Then $C$ lies on a unique smooth quintic del Pezzo surface $\Sigma\subset \PR^5$, embedded by $|-K_{\Sigma}|$, and $C$ is the zero locus of a section of $\OO_{\Sigma}(2)$.  The automorphism group of $\Sigma$ is isomorphic to the symmetric group $S_5$ and the action extends to a linear action of $S_5$ on $\PR^5$ preserving $\Sigma$.

The curve $C$ has exactly five $\grd{1}{4}$'s, and the corresponding five $\grd{2}{6}$'s arise from blowing down one of the five sets of four disjoint lines in $\Sigma$ to obtain an irreducible sextic in $\PR^2$ with four nodes or cusps in linear general position.  Conversely, every plane sextic with four nodes or cusps in linear general position arises in this way, and the canonical series $K_C\cong \OO_C(1)$ is cut out on the plane sextic by the plane cubics passing through the four double points.
\end{prop}

Let $V_d=H^0\left(\Sigma,\OO_\Sigma(d)\right)$ for $d=1,2$.  By Kodaira vanishing and Riemann-Roch for surfaces, we have  \[\dim V_1=\chi(\OO_{\Sigma})+\frac{1}{2}\left((-K_{\Sigma})^2-(-K_{\Sigma})\cdot K_{\Sigma} \right)-h^0\left(\Sigma,2K_\Sigma\right) = 1+5-0=6,\]
since $K_{\Sigma}-(-K_{\Sigma})=2K_{\Sigma}$ is negative,  and similarly 
 \[\dim V_2=\chi(\OO_{\Sigma})+\frac{1}{2}\left((-2K_{\Sigma})^2-(-2K_{\Sigma})\cdot K_{\Sigma} \right)-h^0\left(\Sigma,3K_\Sigma\right) = 1+15-0=16.\]
 These numbers could also have been computed as, respectively, the dimension of the vector space of plane cubics passing through four points in linear general position and the dimension of the vector space of plane sextics double at those four points.
 
\subsection{Stable rationality in the case $0,0,0,0,0,n$}
 
We consider first the case of the vanishing sequence $0,1,2,3,4,5+n$, for $1\leq n\leq 4$.    First, we set \[F=\left\{(\tau,p)\in V_1\times \Sigma\!:~~   \tau(p)=0  \right\}.\]
Then $F\rightarrow \Sigma$ is a rank $5$ vector subbundle of the trivial vector bundle $V_1\times \Sigma\rightarrow \Sigma$ over the surface $\Sigma$.  In particular, $F$ is irreducible of dimension $7$.  Now, for a general $(\tau,p)\in F$ over a given $p\in \Sigma$, the zero locus $T=(\tau)$ of $\tau$ is smooth and irreducible: this follows from Bertini's theorem and the smoothness of $\Sigma$.  Such curves $T$ have genus $1$, as can be checked by adjunction, recalling that $V_1=|\OO_{\Sigma}(1)|= |-K_{\Sigma}|$, so that $g(T) = 1+\tfrac{1}{2}\left((-K_{\Sigma})^2+(-K_{\Sigma})\cdot K_{\Sigma}\right)=1$.  Moreover, by considering the projection $F\rightarrow V_1$, we see that for a general $(\tau,p)\in F$ with $T$ smooth, we have $\OO_T(2)\not\cong \OO_T(10p)$, since for a fixed $T$ there are only $100$ points $p$ on $T$ for which $\OO_T(2)\cong \OO_T(10p)$.

Now, let $F^{\circ}$ be the subset of $F$ consisting of those $(\tau,p)\in F$ for which the zero locus $T=(\tau)$ of $\tau$ is smooth and irreducible with $\OO_T(2)\not\cong \OO_T(10p)$.  Then $F^{\circ}$ is open, since the loci we have excluded are pull-backs of closed subsets of $\PR(V_1)$ and of the Picard scheme over $\M_{1,1}$, and $F^{\circ}\subseteq F$ is dense, since $F$ is irreducible and we have shown $F^{\circ}$ is nonempty.
We set
\[E_n = \left\{(\sigma,\tau,p)\in V_2\times F^{\circ}\!:~~\text{$\sigma$ vanishes to order at least $n+5$ at $p$ along $(\tau)$}\right\},\]
and we claim that $E_n\rightarrow F^{\circ}$ is a vector bundle of rank $16-(n+5)=11-n$.  To show this, we note first that the map $H^0\left(\Sigma,\OO_{\Sigma}\left(2\right)\right)\rightarrow H^0\left(T,\OO_{T}\left(2\right)\right)$ is surjective (this follows from the long exact sequence in cohomology for  $0\to \mathcal{I}_T(2)\to\OO_{\Sigma}(2)\to\OO_{T}(2)\to 0$, since $\mathcal{I}_T(2)\cong\OO_{\Sigma}(2-1)\cong -K_{\Sigma}$).  Now,  $\OO_T(2)$ is a degree $10$ line bundle on $T$, so by Riemann-Roch $h^0(T,\OO_T(2))=10$, and similarly $h^0\left(T,\OO_T\left(2\right)\!\left(-\left(n+5\right)p\right)\right)=10-n-5$, since $n\leq 4$ and $K_T=0$.  Thus vanishing to order $n+5$ at $p$ along $T$ imposes exactly $n+5$ conditions.\footnote{It is here that the $1\leq n\leq 4$ case we have been considering differs from the $n=5$ case:  given a smooth hyperplane section $H$ of $\Sigma$ through a point $p\in\Sigma$, the condition of vanishing along $H$ with multiplicity at least $10$ no longer imposes exactly $10$ independent conditions on $H^0(\Sigma,\OO_{\Sigma}(2))$.  If $C$ is the zero-locus of a global section of $\OO_{\Sigma}(2)$ that does not contain $H$ but meets $H$ with multiplicity at least $10$, then $p$ is the only point of intersection (with multiplicity exactly $10$) between $H$ and $C$ and we must have $\OO_H(2)\cong \OO_H(10p)$.  Since $H$ is a smooth curve of genus $1$, the divisor $q-p$ is not effective for $q\not=p$, so a global section of $\OO_H(10p)$ vanishing to order $9$ at $p$ must in fact vanish to order $10$.  Note that for general $p\in H$, we have $\OO_H(2)\not\cong \OO_H(10p)$, so no smooth section of $H^0(\Sigma,\OO_{\Sigma}(2))$ meets $H$ at $p$ with multiplicity $10$.  For more a more detailed explanation in terms of the corresponding plane cubics and sextics, see Section 4.5 of \cite{mythesis}.}

Thus, $E_n$ is irreducible of dimension $18-n$.  Moreover, by Bertini's theorem, the zero locus $S=(\sigma)$ for a general $(\sigma,\tau,p)\in E_n$ over $(\tau,p)$ is smooth away from $T$, and by considering multiples of $\tau$ by another section of $\OO_{\Sigma}(1)$, we find that a general such $S$ is also smooth along $T$.  By adjunction, such curves $S$ have $g(S)=1+\tfrac{1}{2}\left((-2K_{\Sigma})^2+(-2K_{\Sigma})\cdot K_{\Sigma}\right)=6$ and are canonically embedded in $\PR^5$.  Moreover, a general such smooth $S$ meets $T$ at $p$ to order exactly $n+5$: for $n<4$ this follows since by the dimension count $E_1\supsetneq E_2 \supsetneq E_3\supsetneq E_4$, and for $n=4$ it follows from the condition $\OO_T(2)\not\cong \OO_T(10p)$ on $F^{\circ}$.

Let $E_n^{\circ}$ be the dense open subset of $E_n$ consisting of $(\sigma,\tau,p)$ where $S$ is smooth and irreducible and meets $T$ at $p$ with multiplicity exactly $n+5$.  Then for any element of $E_n^{\circ}$, the vanishing sequence of $K_S$ at $p$ is exactly $0,1,2,3,4,n+5$: we know that $\tau$ vanishes to order exactly $n-5$, and any other section of $\OO_\Sigma(1)$ vanishes on $T$ at five points, counting multiplicities, but cannot vanish at $p$ with multiplicity $5$ because that would imply $\OO_T(1)\cong \OO_T(5p)$ and hence $\OO_T(2)\cong \OO_T(10p)$.

Now, we have a map $E_n^{\circ}\rightarrow \M_{6,1}^H$, where $H$ is the semigroup corresponding to the vanishing sequence $0,1,2,3,4,n+5$, defined by mapping $(\sigma,\tau,p)$ to the class of $(S,p)$.  Since an open subset of $\M_{6,1}$ can be realized by sections of $\OO_{\Sigma}(2)$ (and since the condition that $T$ be smooth is likewise open), we see that the irreducible variety $E_n^{\circ}$ of dimension $18-n$ must dominate a component of $\M_{6,1}^H$.  The fibers of this map are two-dimensional, coming from the $(\C^*)^2$ action by scalar multiplication on $(\sigma,\tau)$, so we find that this component has dimension $16-n=\dim \M_{6,1}-n$, and since $n$ is the weight of these Weierstrass points, this component of $\M_{6,1}^H$ is dimensionally proper.

We now show that this component is stably rational.  Recall that the action of $S_5$ on $\Sigma$ is almost free (\cite{M61}, Claim 1.2).  As in the proof of the rationality of $\M_{6,1}$ in \cite{M61}, our main tool in proving stable rationality will be the following proposition.
\begin{prop}[\cite{Do}, Main Lemma of Section 4]\label{do-mainlemma} Let $G$ be a reductive algebraic group acting almost freely on an irreducible variety $X$ and $E$ be a $G$-linearized vector bundle on $X$.  Then $E/G$ is birationally isomorphic to the total space of a vector bundle on $X/G$.\end{prop}

The action of $S_5$ on $\Sigma$ induces a linearized action on $F$, which is also almost free and preserves $F^{\circ}$, and this in turn induces an $S_5$-linearized action on $E_n$ that preserves $E_n^{\circ}$.  It follows then from Proposition \ref{do-mainlemma} that $E_n^{\circ}/S_5$, which is birational to $E_n/S_5$, is birational to the total space of a vector bundle over $F^{\circ}/S_5$.  Applying Proposition \ref{do-mainlemma} again, we see that $F^{\circ}/S_5$ is birational to the total space of a vector bundle on $\Sigma/S_5$, but a unirational surface is rational, so this implies that $F^{\circ}/S_5$, and hence $E_n^{\circ}/S_5$ is rational.  On the other hand, if $X_n$ is our component of $\M_{6,1}^H$, then $X_n\times \PR^2$ is birational to $E_n^{\circ}/S_5$.  It follows that our component $X_n$ is stably rational.

\begin{rem}Unfortunately, unlike the representation of $S_5$ on $H^0\left(\Sigma,\OO_{\Sigma}(2)\right)$, the representation on
$H^0\left(\Sigma,\OO_{\Sigma}(1)\right)$ does not contain the trivial or sign representation: indeed, it is the $6$-dimensional irreducible representation of $S_5$ (\cite{SB89}, Lemma 1).  We therefore are unable to use the same trick as in \cite{M61} for obtaining rationality of projective bundles from rationality of vector bundles.
\end{rem}

\subsection{Stable rationality in the case $0,0,0,0,1,n$}

We now consider the case where $H$ is the semigroup corresponding to the vanishing sequence $0,1,2,3,5,5+n$, for $n=1,2$.  If $S$ is the smooth irreducible zero-locus of a section $\sigma$ of $\OO_{\Sigma}(2)$, and $p$ is a Weierstrass point of $S$ with semigroup $H$, then
there will exist sections $\tau$ and $\upsilon$ of $\OO_{\Sigma}\left(1\right)$ vanishing to orders $5+n$ and $5$, respectively, along $S$ at $p$.  While $\tau$ is uniquely determined up to scalars as in the previous case, $\upsilon$ is only uniquely determined up to scalars modulo $\tau$.  We note, however, that as long as $T=(\tau)$ is smooth, there is a unique (up to scalars) linear combination $\upsilon+c\tau$ whose zero locus $U$ is singular at $p$; we will choose that as our $\upsilon$.

The quintic del Pezzo surface $\Sigma$ contains ten lines; if we regard $\Sigma$ as the blow-up of $\PR^2$ at four points, these are the proper transforms of the six lines passing through pairs of these points, together with the four exceptional divisors.  Let $\Sigma^{\circ}$ be the complement of these ten lines in $\Sigma$.

  Let $A=\{(w,p)\!:p\in \Sigma^{\circ}, w\in T_p\Sigma^{\circ}\} \rightarrow \Sigma^{\circ}$ be the tangent bundle of $\Sigma^{\circ}$, and let $A^{\circ}$ be the open subset where $w\not=0$.  Then $A^\circ$ is irreducible of dimension $4$.  Now, set
\[B=\left\{\left(\upsilon,w,p\right)\in V_1\times A^{\circ}\!: ~\mult_p(\upsilon)\geq 2\text{ and $\upsilon$ vanishes to order at least $3$ in the direction $w$}  \right\},\]
where the condition of vanishing to order at least $3$ in the direction given by a vector is well-defined because we are assuming that $\upsilon$ vanishes to order at least $2$ in all directions; in the case where the zero locus $U=(\upsilon)$ of $\upsilon$ has a double point at $p$, we are simply requiring that some branch have tangent direction $w$.  Then $B\rightarrow A^{\circ}$ is a rank $6-4=2$ vector bundle\footnote{One can check directly that given five distinct points in $\PR^2$ and a tangent vector at the fifth point, the $4+3+1$ conditions on the space of plane cubics of vanishing at the first four points and vanishing with multiplicity at least $2$ in all directions and at least $3$ in the specified direction at the fifth point are independent as long as the five points are in linear general position.} and hence $B$ is irreducible of dimension $6$.

Let $B^{\circ}$ be the dense open subset of $B$ where $U=(\upsilon)$ is irreducible with a node or cusp at $p$.  Set
\[C=\left\{\left(\tau,\upsilon,w,p\right)\in V_1\times B^{\circ}\!:~~\text{$(T.U)_p\geq 5$ and $\tau$ vanishes to order at least $2$ at $p$ along $w$}  \right\},\]
where $T=(\tau)$ is the zero-locus of $\tau$, and $(T.U)_p$ is the intersection multiplicity at $p$, but considered to be $+\infty$ if $\tau$ is identically zero on $U$.  Then $C\rightarrow B^{\circ}$ is a vector bundle of rank $2$.\footnote{Let $\fun{\pi}{\PR^1}{U}$ be the normalization of the nodal or cuspidal rational curve $U$.  Then $\pi^*\OO_U(1)\cong \OO_{\PR^1}(5)$, and $\pi^*$ induces a linear map $V_1\rightarrow H^0\left(\PR^1,\OO_{\PR^1}\left(5\right)\right)$ between $6$-dimensional vector spaces whose kernel is the $1$-dimensional subspace spanned by $\upsilon$.  Thus the image is precisely the $5$-dimensional subspace defined by the $\delta=1$ adjoint condition (see \cite{ACGH} I.A.2) and in both the node case and the cusp case we find that we are imposing exactly $4$ linear conditions on the $6$-dimensional vector space $V_1$.  (In the node case, vanishing to order $4$ along the specified branch guarantees an intersection multiplicity of at least $5$, and in the cusp case, the possible vanishing orders at the point on $\PR^1$ corresponding to $p$ are $0,2,3,4,5$.) 
}
Therefore $C$ is irreducible of dimension $8$.  Let $C^{\circ}$ be the dense open subset where $T$ is smooth, and set
\[D_n=\left\{\left(\sigma,\tau,\upsilon,w,p\right)\in V_2\times C^{\circ}\!:~~\text{$\sigma$ vanishes to order $5+n$ along $T$ at $p$}  \right\} \]
for $n=1,2,3$.  As in the case of $E_n\rightarrow F^{\circ}$ above, we find that $D_n\rightarrow C^{\circ}$ is a rank $11-n$ vector bundle, and for a general $(\sigma,\tau,\upsilon,w,p)\in D_n$, the curve $S=(\sigma)$ is a smooth canonically embedded curve of genus $6$.  We have then that $D_n$ is irreducible of dimension $19-n$, and since this implies $D_1\supsetneq D_2 \supsetneq D_3$, we find that at a general point of $D_n$ for $n=1,2$, the intersection multiplicity of $S$ and $T$ at $p$ is exactly $5+n$.

Let $D_n^{\circ}$ be the dense open subset where $S$ is a smooth irreducible curve and $S$ and $T$ meet with multiplicity exactly $5+n$ at $p$.  Then since $T$ and $U$ meet at $p$ with multiplicity exactly $5$, we must have that $S$ and $U$ meet at $p$ with multiplicity exactly $5$.  Moreover, since $C\rightarrow B^{\circ}$ has rank $2$, the sections $\sigma$ and $\upsilon$ are a basis  for  the subspace of $V_1$ vanishing to order at least $5$ along $S$.  On the other hand, no section of $\OO_{\Sigma}(1)$ can vanish to order exactly $4$ along $S$ at $p$, since then on the smooth curve $T$ of genus $1$ we would have $\OO_T(4p+q)\cong \OO_T(1) \cong \OO_T(5p)$ for some point $q\neq p$ of $T$. It follows that for every $(\sigma,\tau,\upsilon,w,p)\in D_n^{\circ}$, with $n=1,2$, the point $p$ is a Weierstrass point of $S=(\sigma)$ with vanishing sequence $0,1,2,3,5,5+n$.

We thus have, for $n=1,2$, a map $D_n^{\circ}\rightarrow\M_{6,1}^H$, where $H$ is the semigroup corresponding to the vanishing sequence $0,1,2,3,5,5+n$.  As in the previous case, we find that this map dominates some component of $\M_{6,1}^H$.  The fibers of this map are four-dimensional, coming from the $(\C^*)^4$ action by scalar multiplication on $(\sigma,\tau,\upsilon,w)$, so this component has dimension $15-n=16-(n+1)$, and the component is dimensionally proper.  The proof of the stable rationality of this component $Y_n$ then proceeds exactly as in the previous case.  We apply Proposition \ref{do-mainlemma} four times to conclude that $D_n^{\circ}/S_5$ is rational, and it follows that $Y_n\times \PR^4$ is rational.

\begin{rem} This entire argument also goes through verbatim when $n=3$ for the vanishing sequence $0,1,2,3,5,8$, but in this case the corresponding semigroup is $N(6)_{13}=\langle 5,7,8,11\rangle$, and the irreducibility and rationality of $\M_{6,1}^{N(6)_{13}}$ was proved in \cite{Nak2}.
\end{rem}

\subsection{Irreducibility}

At this point, we have shown only the existence of an irreducible component of the expected dimension for each semigroup; it's conceivable there could be another component that simply fails to meet the open subset that we've been considering, e.g. because it consists only of trigonal curves, bielliptic curves, or the other excluded cases in Proposition \ref{acghprop}, or because for a general $(C,p)$ in the component, the zero locus $T=(\tau)$ of the section of $\OO_{\Sigma}(1)$ meeting $C$ with multiplicity $5+n$ is singular rather than smooth.  However, since every component must have codimension at most the weight of the semigroup, to show that we haven't missed a component we need only repeat the same sort of dimension count (but without worrying about rationality) in each special case and show that we get a dimension smaller than $16-n$ in the first case and $16-(n+1)$ in the second.

The Weierstrass points in the hyperelliptic case have vanishing sequence $0,2,4,6,8,10$, which is not one of the sequences we are considering.  The trigonal case has also been worked out completely (in arbitrary genus) in \cite{trigonal}, Theorem 4.6.\footnote{The semigroups we are concerned with can all occur in the trigonal case.
In the notation of \cite{trigonal}, the vanishing sequence $0,1,2,3,5,7$ corresponds to $\rho=s_E=m=1$, $t_E=t=\epsilon=0$, yielding a dimension of $11<13$.  The vanishing sequence $0,1,2,3,5,6$ corresponds to $\rho=\epsilon=\alpha=1$, $t=t_E=s_E=0$, $m=2$, $r=4$, yielding a dimension of $12<14$.  For $n=2,3,4$, the vanishing sequence $0,1,2,3,4,5+n$ corresponds to $\rho=\epsilon=1$, $t=t_E=s_E=0$, $m=2$, $r=5$, $\alpha=n$, yielding a dimension of $14-n<16-n$.}  The entire locus $\M_6^{\text{be}}$ of bielliptic curves of genus $6$ has dimension only $2\cdot6-2=10$, so it is certainly too small to contain a component for any of the semigroups we are considering.  

The locus of smooth plane quintics (up to projective equivalence) has dimension only $\binom{5+2}{2}-1-\dim PGL_3=12$, so we need only check that it is not the case that a general smooth plane quintic has a Weierstrass point of vanishing sequence $0,1,2,3,4,9$.  But let $C\subset \PR^2$ be a smooth quintic and $p\in C$ be a Weierstrass point with vanishing sequence $0,1,2,3,4,9$.  Then (by adjunction) the canonical series $|K_C|$ is cut out on $C$ by the plane conics.  If $Z$ is the conic meeting $C$ with multiplicity $9$ at $p$, then $Z$ can't be reducible, since then the double tangent line would meet $C$ with a greater multiplicity, and $Z$ can't be a double line since $9$ is odd.  Thus $Z$ is smooth, and we can easily count the dimension of the the image in $\M_{6,1}$ of the locus of pointed curves $(C,p)$ of this form:
\[\overbrace{2}^{\text{point $p$}} + \overbrace{\binom{2+2}{2}-1}^{\text{smooth conic $Z$,}} ~- \overbrace{1}^{\text{through $p$}} + \overbrace{\binom{5+2}{2}-1}^{\text{smooth quintic $C$,}}-\overbrace{9}^{\text{with $(C.Z)_p = 9$}} - \overbrace{8}^{\text{$PGL_3$}}-\overbrace{0}^{\dim G^2_5(C)} = 9,\]
which is less than $16-4=12$, the smallest possible dimension of a component in this case.

The remaining special cases in Proposition \ref{acghprop} are those where a general curve $S$ in a possible component is birational to a singular sextic curve $\widetilde{S}\subseteq \PR^2$ with double point singularities, at least one of which is of analytic type $A_m: y^2=x^{m+1}$ for $3\leq m \leq 8$ (cf. \cite{ACGH} V.A).  Since the arithmetic genus of $\widetilde{S}$ is $\binom{6-1}{2}=10$, the $\delta$-invariants of these singularities must add up to $10-6=4$, so we must have either a single singularity with $\delta=\left\lfloor\frac{m}{2}\right\rfloor=2,3,4$ and $4-\delta$ other nodes or cusps, or two singularities with $\delta=2$ (i.e. tacnodes or rhamphoid cusps).  We can handle these cases the same way we handled the general case: we could have expressed our dimension count in terms of sextic plane curves with four nodes or cusps in linear general position, but it was more convenient to instead consider their proper transform in the blowup of $\PR^2$ at those points.  Likewise, in these cases we instead consider the rational surface $\Sigma$ obtained in resolving the singularities of $\widetilde{S}$.

In each case, the points that we blow up at each stage in resolving the singularities of $\widetilde{S}$ do not lie on any irreducible $(-2)$-curves\footnote{If the tangent line to a tacnode were also tangent to a branch of another node (or tacnode), this condition would fail: the proper transform of the tangent line in $\Sigma$ would have self-intersection $-3$.  This case, however, never occurs for a plane sextic by B\'{e}Žzout's theorem.}, so the resulting surface $\Sigma$ that we obtain is a generalized (or ``weak'') Del Pezzo surface of degree $5$ in the sense of \cite{DPfr} or \cite{DPen}.  Likewise, one can check in each case that the proper transform $S$ of $\widetilde{S}$ is the zero locus of a section of $-2K_\Sigma$, with $K_S$ cut out by $-K_\Sigma$.

\begin{figure}
\includegraphics{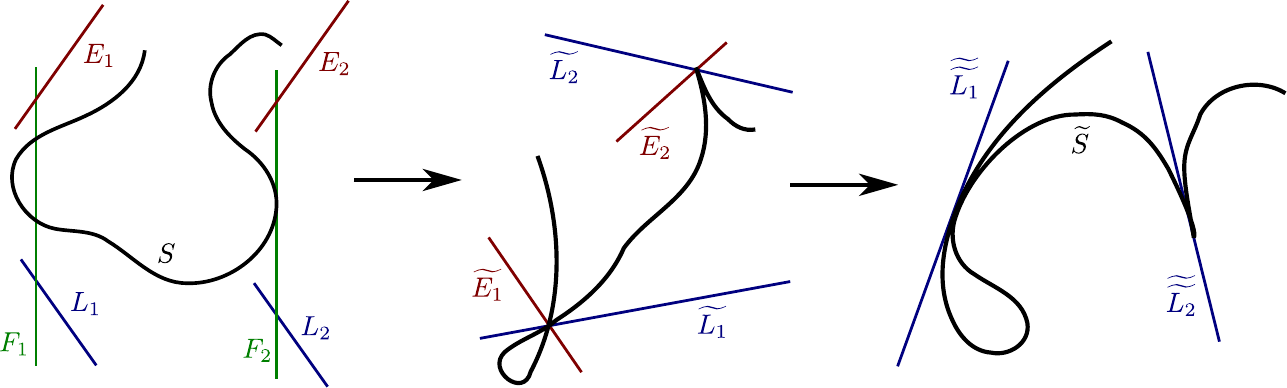}
\caption{The case where the plane sextic curve $\widetilde{S}$ has an $A_3$ and an $A_4$ singularity.}
\end{figure}

For example, in the case of two tacnodes, each singularity is resolved by blowing up at the singular point and then blowing up at the infinitely near point corresponding to the shared tangent line of the branches.  In this case, $\Z^5\cong\Pic(\Sigma)=\langle H, E_1,E_2, F_1, F_2\rangle$, where $H$ is the pull-back of the class of a line in $\PR^2$, $E_1$ and $E_2$ are the proper transforms of the original exceptional divisors, and $F_1$ and $F_2$ are the exceptional divisors from blowing up the infinitely near points.  The $E_i$ are $(-2)$-curves, the $F_i$ are $(-1)$-curves and $E_i\cdot F_j = \delta_{ij}$.  The canonical class on $\Sigma$ is \[K_{\Sigma}=\OO_{\Sigma}(-3H+E_1+2F_1+E_2+2F_2),\]
and he tacnodal plane sextic $\widetilde{S}$ pulls back to a section of $\OO_{\Sigma}(6H)$ vanishing on $2E_1+4F_1+2E_2+4F_2$, so its proper transform $S$ has class $6H-2E_1-4F_1-2E_2-4F_2 = -2K_{\Sigma}$.  Note that $\widetilde{S}$ meets $F_i$ in $-2E_i\cdot F_i-4F_i\cdot F_i=-2+4=2$ distinct points.  In the case that one or both of the singularities were $A_4$ instead of $A_3$, the surface $\Sigma$ and the class of $S$ would be the same, but $S$ would meet the corresponding $F_i$ in one point with multiplicity two rather than two distinct points.

Thus we can repeat the same dimension count as in the case of an ordinary Del Pezzo surface above, replacing $\OO_\Sigma(1)$ and $\OO_\Sigma(2)$  with $-K_\Sigma$ and $-2K_{\Sigma}$.  The line bundle $-K_\Sigma$ is no longer ample, but it is big and nef, so we may use the Kawamat-Viehweg vanishing theorem (\cite{Kvanish}\cite{Vvanish}) where we used Kodaira vanishing above.  The only difference comes at the end when we consider the fibers of the map from our parameter space to $\M_{6,1}$: while the automorphism group of an ordinary Del Pezzo surface of degree $5$ is finite, the automorphism groups of these generalized Del Pezzo surfaces of degree $5$ are all positive-dimensional,\footnote{This is because demanding that an element of $PGL_3$ that is known to fix a given point also fix an infinitely near point only imposes at most one additional condition rather than the two conditions imposed by fixing another general point in $\PR^2$.  If it imposes zero conditions instead of one, then the (coarse) moduli space of such surfaces has larger dimension, but the fact that the automorphism group is even bigger in this case makes up for this in the dimension count.}
 so the dimension of the fibers of the map to $\M_{6,1}$ is larger and the image in $\M_{6,1}$ has smaller dimension.

We must now deal with the cases where the auxiliary curve $T$ is singular.  We can analyze these cases in, for the most part, the same way we analyzed the general case (but without having to show that the general vanishing sequence isn't larger than we expect or even that the general curves in our families are actually smooth).  We discuss briefly here the modifications required.  First of all, if $T$ is irreducible and singular, but not singular at $p$, then the analysis goes through unchanged (but the proof that vanishing to order $5+n$ at $p$ imposes $5+n$ independent conditions involves a smooth point of a nodal or cuspidal rational curve of arithmetic genus $1$ rather than a smooth curve of geometric genus $1$) producing a locus of one lower dimension.   We can handle the cases where $T$ is reducible but not singular at $p$ in the same way.

We now consider the cases where $T$ is singular at $p$.  When $T$ has a node or cusp at $p$, the requirement that $(S.T)_p\geq 5+n$ imposes only $4+n$ conditions,\footnote{In the node case, the fibers of this map are no longer irreducible, but we can fix this by choosing a tangent direction for one of the branches of the node first.} but requiring a singularity at $p$ reduces the dimension of the possible choices of $(T,p)$ by two.\footnote{For ramification sequences $0,0,0,0,1,n$, a further slight modification is required: rather than picking $U$ first and then $T$, we pick $T$ first and then $U$; it's no longer the case that $U$ us uniquely determined by $T$, but since we're just proving an upper bound on the dimension, we don't need to keep track of this.}  When $T$ has a tacnode (this requires that $T$ be reducible), the requirement that $(S.T)_p\geq 5+n$ imposes only $3+n$ conditions, but the dimension of the space of tacnodal $(T,p)$ is smaller by four.\footnote{In fact, the tacnodal case is the general situation for the semigroup $N(6)_5$ with vanishing sequence $0,1,2,4,5,6$.  Identifying $\Sigma$ with the blowup of $\PR^2$ at $q_1,\ldots,q_4$, we consider the plane conic $Z$ through $p$ and the $q_i$, and demand that that the image of $S$ meet $Z$ with multiplicity $4$.  Then the image of $T$ in $\PR^2$ is $Z\cup L$, where $L$ is the tangent line to $Z$ at $p$.}  When $T$ has an ordinary triple point (this requires that $T$ be reducible and that $p$ lie on a line of $\Sigma$), the requirement that $(S.T)_p\geq 5+n$ imposes only $3+n$ conditions, but the dimension of the space of such $(T,p)$ is smaller by $5$.  In the case where $T$ is non-reduced at $p$ but $T_{\text{red}}$ is smooth at $p$,  the requirement that $(S.T)_p\geq 5+n$ imposes only $\left\lceil\frac{5+n}{2}\right\rceil$ conditions, but the dimension of the space of such $(T,p)$ is smaller by $5$ (and $n\leq 4$).  Finally, in the case where $T$ is non-reduced at $p$ and $T_{\text{red}}$ is nodal at $p$, the requirement that $(S.T)_p\geq 5+n$ imposes only $\left\lceil\frac{4+n}{2}\right\rceil$ conditions, but the dimension of the space of such $(T,p)$ is lower by $6$ (indeed, there are only finitely many such $(T,p)$).\footnote{In the case of the ramification sequences $0, 0, 0, 0, 1, n$, we must also deal with the case where $p$ is on a line of $\Sigma$ and the vector $w$ points along the line, but in fact this implies $T$ is non-reduced.  Any section of $\OO_{\Sigma}(1)$ vanishing at $p$ to order at least $2$ in the direction $w$ must vanish on the entire line, meaning that vanishing to order at least $3$ at $p$ in the direction $w$ imposes only $2$ conditions on $U$ instead of $3$, and meeting $U$ with multiplicity $5$ imposes only $2$ conditions on $T$ instead of $4$.  On the other hand, we can't just choose $T$ arbitrarily according to the same conditions as usual: given that $T$ must contain the line and must have the biggest possible intersection multiplicity with $S$ at $p$, in fact $T$ must contain the double line, and we are in a non-reduced case considered above.}  This completes the proof of the theorem.

\begin{rem} There is no real reason to believe that the two remaining cases are not stably rational (or, indeed, that any of these varieties fail to be rational).  Indeed, one might expect the behavior in the case of the semigroup $\langle 6,7,8,9,10\rangle$ of a general odd subcanonical point (\cite{mythesis}) in genus $6$ to closely resemble the case of the semigroup $\langle 5,7,8,9\rangle$ of a general even subcanonical point, for which rationality is proved in \cite{Nak2}.

\end{rem}

\end{document}